\documentclass[11pt]{article}

\usepackage{multirow, multicol}
\usepackage{booktabs} % for much better looking tables
\usepackage{array} % for better arrays (eg matrices) in maths
\usepackage{paralist} % very flexible & customisable lists (eg. enumerate/itemize, etc.)
\usepackage{verbatim} % adds environment for commenting out blocks of text & for better verbatim
\usepackage{subfig} % make it possible to include more than one captioned figure/table in a single float
\usepackage{amsmath, amssymb, enumerate, cite}
\usepackage{pictexwd,dcpic}

\usepackage[final]{graphicx}

\newtheorem{theorem}{Theorem}

\newtheorem{remark}[theorem]{Remark}
\newtheorem{lemma}[theorem]{Lemma}
\newtheorem{definition}[theorem]{Definition}
\newtheorem{proposition}[theorem]{Proposition}
\newtheorem{example}[theorem]{Example}

\newcommand{\proof}{ {\sc Proof.\quad}}
\newcommand{\pend}{ \hfill $\square$ \\}

\numberwithin{equation}{section}  % Formeln mit f¸hrender Kapitelnummer
\numberwithin{figure}{section}    % Abbildungen mit f¸hrender Kapitelnummer
\numberwithin{table}{section}     % Tabellen mit f¸hrender Kapitelnummer
\numberwithin{theorem}{section}

\newcommand{\of}[1]{\ensuremath{\left( #1 \right)}}

\newcommand{\cb}[1]{\ensuremath{ \left\{ #1 \right\} }}
\newcommand{\sqb}[1]{\ensuremath{ \left[ #1 \right] }}

\newcommand{\st}{\,|\;}

\newcommand{\eps}{\ensuremath{\varepsilon}}

\newcommand{\vp}{\ensuremath{\varphi}}

\newcommand{\R}{\mathrm{I\negthinspace R}}
\newcommand{\OLR}{\overline{\mathrm{I\negthinspace R}}}
\newcommand{\N}{\mathrm{I\negthinspace N}}

\renewcommand{\P}{\ensuremath{\mathcal{P}}}

\newcommand{\G}{\ensuremath{\mathcal{G}}}

\newcommand{\C}{\ensuremath{C^-\setminus\cb{0}}}

\newcommand{\dom}{{\rm dom \,}}

\newcommand{\gr}{{\rm graph \,}}
\newcommand{\cl}{{\rm cl \,}}

\newcommand{\co}{{\rm co \,}}

\newcommand{\cone}{{\rm cone\,}}

\newcommand{\Int}{{\rm int\,}}

\newcommand{\isum}{{+^{\negmedspace\centerdot\,}}}

\newcommand{\idif}{{-^{\negmedspace\centerdot\,}}}

\newcommand{\lel}{\preccurlyeq}

\newcommand{\triup}{{\rm \vartriangle}}

\usepackage[%backref=page,
colorlinks=true, linkcolor=blue, citecolor=blue, urlcolor=blue, filecolor=blue
]{hyperref}

\title{Weak Minimizers, Minimizers and Variational Inequalities for set--valued Functions. A blooming wreath?}

\author{
Giovanni P. Crespi\thanks{University of Insubria, Department of Economics,  Via Monte Generoso, 71, 21100, Varese, Italy, {giovanni.crespi@uninsubria.it}}
and
Carola Schrage\thanks{Free University of Bozen-Bolzano,
Faculty of Economics and Management, Piazza Universita 1, 39100 Bozen-Bolzano Italy, {carolaschrage@gmail.com}}%\footnote{the research was partly done while the second authors stay at the University of Valle d'Aosta}
}
\begin{document}
\maketitle
\begin{abstract}
In \cite{CrespiSchrage13a}, necessary and sufficient conditions in terms of variational inequalities are introduced to characterize minimizers of convex set--valued functions with values in a conlinear space. Similar results are proved in \cite{CrespiSchrage13b,CrespiRoccaSchrage14W} for a weaker concept of minimizers and weaker variational inequalities. The implications are proved using scalarization techniques that eventually provide original problems, not fully equivalent to the set-valued counterparts.
Therefore, we try, in the course of this note, to close the network among the various notions proposed. More specifically, we prove that a minimizer is always a weak minimizer, and a solution to the stronger variational inequality always also a solution to the weak variational inequality of the same type. As a special case we obtain a complete characterization of efficiency and weak efficiency in vector optimization by set-valued variational inequalities and their scalarizations. Indeed this might eventually prove the usefulness of the set optimization approach to renew the study of vector optimization.


\medskip{\bf Classcode:}
49J40, 	49J53 , 	58C06, 58E30
\end{abstract}

\section{Introduction}\label{sec:Intro}

Scalar variational inequalities (for short VI) apply to study a wide range of problems, such as equilibrium and optimization problems,  see e.g. \cite{BaiocchiCapelo84, KinderlehrerStempacchia80}.
Generalizations toward vector VI were initiated in \cite{Giannessi80}; for recent results and surveys on this field see  \cite{Giannessi98, GiannessiMastroeniPellegrini2000, Komlosi99, LiMastroeni2010}. A major peculiarity in vector--valued inequalities is the necessity to introduce at least two different solution concepts, e.g. a strong and a weak one. This approach seems to be most natural if referred to vector optimization efficiency and weak efficiency notions.

%Far less has been undertaken to extend those results to set-valued optimization, mainly because of a rather different approach to the classical optimization of set-valued maps.

The notion of differentiable variational inequalities arises in the scalar case, when the operator involved in a VI has a primitive function. This kind of VI is widely studied because of its relation to optimization problems. Under mild continuity assumptions, scalar Minty VI (MVI, \cite{LionsStampacchia67, Minty67}) of differential type provide a sufficient optimality condition to the primitive optimization problem (a result popularized as Minty Variational Principle), while scalar Stampacchia VI (SVI, \cite{Stampacchia64}) is only necessary. Assuming some convexity on the primitive function (or monotonicity of the derivative) both VIs are necessary and sufficient optimality conditions. In \cite{crespiJOTA}, under generalized differentiability assumptions, scalar Minty VI have been studied and it has been proved that the existence of a solution to such a problem implies some regularity property on the primitive optimization problem.

The same approach has been proposed by Giannessi \cite{Giannessi80} for vector optimization. In his seminal paper, Giannessi studied the relations between a Stampacchia type vector variational inequality and weak efficient solutions of the primitive vector optimization problem. It has been proved that the scalar relations hold under stronger assumption in the vector case, namely convexity plays a bigger role in the proof. Further researches tried to extend the result to efficient solutions, providing a stronger version of the variational inequality and highlighting some peculiarities of the vector case unknown for scalar functions.
More recently, also the Minty Variational Principle has been studied and extended to the vector case. The problem has been posed by Giannessi in \cite{Giannessi98}, where the links between Minty VIs and vector optimization problems were investigated both for efficient and weak efficient solutions. More recently, in \cite{CreGinRoc08, YangYang04},  some generalization of the vector variational principle have been proposed in conjuction with weak efficient solutions. In \cite{Giannessi98, YangYang04}, the case of a differentiable objective function $f$ with values in ${\mathbb R}^m$ and a Pareto ordering cone has been studied, proving a vector Minty variational principle for pseudoconvex functions. In \cite{CreGinRoc08} a similar result has been extended to the case of an arbitrary ordering cone and a non-differentiable objective function. Overall, the existing literature pictures a wreath of relations, ranging between weak and strong vector--valued inequalities and weak efficiency and efficiency. Some of these relations occur only under (generalized) convexity assumptions, some of the branches cannot be fixed.

Although optimization of set-valued functions has been a fast growing topic over the past decades, very few has been proposed about variational inequalities to characterize minimality.

Since the first results by Corley \cite{corley87, Corley88} and Dinh The Luc \cite{Luc89}, based on a vector optimization approach, several papers have been proposed to provide optimality conditions. Nevertheless, the main approach to derivatives (and therefore to the core of a variational inequality) has been far distant form the basic differential quotient method adopted for scalar (and vector) problems. More recently, a new paradigm, known as set optimization, has been proposed, compare \cite{HamelHabil05, HeydeLoehne11, LoehneDiss, Loehne11Book}. In this framework, the very concept of optimal solutions has been thought anew, together with operations among sets, now elements of a complete ordered conlinear space. This leads to overcome some drawbacks in previous attempts to provide variational inequalities for set-valued optimization problems (see e.g. \cite{CreGinRoc2010}).

%In the present approach, we consider two types of variational inequalities. First, a set of scalarizations of a set--valued function is given and the variational inequality is expressed using the (extended) Dini derivative of these scalarizations. Secondly, a set--valued Dini derivative of the function itself is defined and used in the variational inequalities.

In \cite{CrespiRoccaSchrage14W} and \cite{CrespiSchrage13b}, a notion of weak minimality for set optimization is presented, motivated by its relation with standard weak efficiency in vector optimization. Under certain regularity assumptions it is proven in \cite{CrespiRoccaSchrage14W} that the solutions of the Minty type inequality are weak minimizers of the primitive set optimization problem. Under slightly weaker assumptions, a weak minimizer of the set optimization problem solves the Stampacchia differential variational inequality. Under convexity assumptions on the scalarizations, the reverse implications has been proven in \cite{CrespiSchrage13b}.
In \cite{CrespiHamelSchrage13W} and \cite{CrespiSchrage13a}, a corresponding chain of implications has been provided for minimizers, actually for solutions of set optimization problems, and the corresponding Minty and Stampacchia type differential variational inequalities.

The aim of this paper is to combine lose branches from the previous studies  between set optimization and set-valued variational inequalities, connecting strong notions in \cite{CrespiSchrage13a} with their weak counterparts presented in \cite{CrespiSchrage13b}. As a special case of our results, we obtain a wreath containing vector optimization efficient and weak efficient solutions.

The paper is organized as follows. We present the general setting of the problem and the basic notation and assumption in Section 2, where some details on conlinear spaces are recalled. In Section 3 we introduce the notion of minimizer and weak minimizer in set optimization, as well as the scalarization technique that is used to prove the main results. The variational inequalities introduced in \cite{CrespiSchrage13a} and \cite{CrespiSchrage13b} are also recalled together with the chains of implications proved in these  papers. Section 4 completes the picture with the main results proving the missing implications.

When of interest, counterexamples are included to show that assumptions cannot be relaxed. The chains of implication provided in each section are illustrated by diagrams.

%%%%%%%%%%%%%%%%%%%%%%%%%%%%%%%%%%%%%%%%%%%%%%%%%%%%%%
\section{Basics}

Throughout the paper, $X$ and $Z$ are locally convex and Hausdorff with topological duals $X^*$ and $Z^*$. The set $\mathcal U_X$, $\mathcal U_Z$ are the set of all closed, convex and balanced $0$--neighbourhoods, respectively, in $X$ and $Z$, that is $\mathcal U_X$ and $\mathcal U_Z$ are $0$--neighbourhood bases. By $\cl A$, $\co A$ and $\Int A$, we denote the closed or convex hull of a set $A\subseteq Z$ and the topological interior of $A$, respectively. The conical hull of a set $A$ is $\cone A=\cb{ta\st a\in A,\, 0<t}$. A cone is Daniell if any bounded decreasing net converges to its infimum.

 The set $Z$ is preordered by a closed convex cone $C\neq Z$ with nonempty topological interior, $\Int C\neq\emptyset$, by means of $z_1\leq_C z_2$ if $z_2\in \cb{z_1}+C$. The (negative) dual cone of $C$ is the set $C^-=\cb{z^*\in Z^*\st \forall z\in C:\, z^*(z)\leq 0}$. Since $\Int C\neq \emptyset$, there exists a weak$^*$ compact base $B^*$ of $C^-$, i.e. a convex subset with $\C=\cone B^*$ with $z^*, tz^*\in B^*$ implying $t=1$ and any net in $B^*$ has a weak$^*$ convergent subnet, compare \cite[Theorem 1.5.1]{Aubin71}

In the sequel we consider the family of subsets of $Z$
\[
\G(Z,C)=\cb{A\in\P(Z)| A=\cl\co(A+C)}.
\]
According to the order relation
\[
\forall A,B\in \G(Z,C):\quad \of{A\lel B \; \mbox{iff} \; B\subseteq A } 
\]
the set $(\G(Z,C),\lel)$ is order complete. %, having set  iff $B\subseteq A$ for all $A,B\in \G(Z,C)$.
Indeed, for any subset $\mathcal A\subseteq \G(Z,C)$ it holds
\[
\inf\mathcal A=\cl\co\bigcup\limits_{A\in\mathcal A}A;\quad \sup\mathcal A=\bigcap\limits_{A\in \mathcal A}A.
\]
Assuming, by definition, that when $\mathcal A=\emptyset$ we have $\inf\mathcal A=\emptyset$ and $\sup\mathcal A =Z$.
Particularly, $\G(Z,C)$ possesses a least element $\inf\G(Z,C)=Z$ and a greatest one $\sup\G(Z,C)=\emptyset$.

We can also introduce algebraic operations on $\G(Z,C)$ by
\begin{align*}
\forall A, B\in \G(Z,C):\quad 		&A\oplus B=\cl\cb{a+b\in Z\st a \in A,\, b\in B};\\
\forall A\in\G(Z,C),\, \forall 0<t:\quad &t\cdot A=\cb{ta \in Z\st a\in A };\quad 0\cdot A=C,
\end{align*}
The resulting space $\G^\triup=(\G(Z,C),\oplus,\cdot,C,\lel)$ is endowed with neutral element $C$, $\emptyset$ dominates the addition and $0\cdot\emptyset=0\cdot Z=C$. Moreover,
\begin{align*}
\forall \mathcal A\subseteq \G(Z,C),\, \forall B\in\G(Z,C) :\quad B\oplus \inf\mathcal A=\inf\cb{B\oplus A\st A\in \mathcal A},
\end{align*}
or, equivalently, the $\inf$--residual 
$$A\idif B=\inf\cb{M\in\G(Z,C)\st A\lel B\oplus M}$$
 exists for all $A, B\in \G(Z,C)$. It holds (compare \cite[Theorem 2.1]{HamelSchrage12})
\begin{align*}
A\idif B    &=\cb{z\in Z\st B+\cb{z}\subseteq A};\\
A             &\lel B\oplus (A\idif B).
\end{align*}

Overall, the structure of $\G^\triup$ is that of an order complete $\inf$--residuated conlinear space as given in the following definition, compare also \cite{HamelSchrage12}.
\begin{definition}
\label{DefConlinearSpace} A nonempty set $Y$ together with two algebraic operations $\isum
: Y \times Y \to Y$ and $\cdot : \R_+ \times Y \to Y$ is called a conlinear
space with neutral element $\theta$ provided that
\\
(C1) $\of{Y, \isum,\theta}$ is a commutative monoid with neutral element $\theta$,
\\
(C2) The operations are compatible: (i) $\forall y_1, y_2 \in Y$, $\forall r \in \R_+$: $r \cdot \of{y_1 \isum y_2} = r
\cdot y_1 \isum r \cdot y_2$, (ii) $\forall y \in Y$, $\forall r, s \in \R_+$: $s \cdot \of{r
\cdot y} = \of{rs} \cdot y$, (iii) $\forall y \in Y$: $1 \cdot y = y$, (iv)  $\forall y \in Y$: $0 \cdot
y = \theta$.

A conlinear space $\of{Y, \isum, \cdot,\theta}$ together with an order relation $\lel$ on $Y$ is called partially ordered, lattice ordered or order complete  conlinear space provided
that $(Y,\lel)$ has the respective structure and  the order is compatible with the algebraic operations $\isum$ and $\cdot$, that is\\
(C3) (i)$\forall y, y_1, y_2 \in Y$, $y_1 \lel y_2$ implies $y_1 \isum y \lel y_2 \isum y$,
 and (ii) $\forall y_1, y_2 \in Y,\,\forall r\in\R_+$, $y_1 \lel y_2$ implies $r \cdot y_1 \lel r\cdot y_2$.

A partially ordered conlinear space  $\of{Y, \isum, \cdot,\theta,\lel }$ is called $\inf$--residuated, when for all $v,y\in Y$ the element
$y\idif v=\inf\cb{u\in Y\st y\lel v\isum u}$
exists. In this case, $y\idif v$ is called the $\inf$--residual of $y$ and $v$.
\end{definition}

We refer to \cite{Fuchs66, galatos2007residuated, GetanMaLeSi, HamelSchrage12, HamelHabil05, MartinezLegazSinger95} for a more thorough study of this structure. For the sake of completeness, we recall that it can be proven that a partially ordered conlinear space is $\inf$--residuated if and only if for all $y\in Y$ and all $A\subseteq Y$ such that $\inf  A$ exists, it holds $\of{y\isum \inf\mathcal A}=\inf\cb{y\isum a\st a\in A}$, compare \cite[Theorem 2.1]{HamelSchrage12}. The structure of a conlinear space is illustrated in  the following example.

\begin{example} \label{ExExtReals}
Let us consider $Z = \R$, $C = \R_+$. Then $\G\of{Z, C} = \cb{[r, +\infty) \mid r \in \R}\cup\cb{\R}\cup\cb{\emptyset}$, and $\G^\triup$ can be identified (with respect to the algebraic and order structures which turn $\G\of{\R, \R_+}$ into an ordered conlinear space and a  complete lattice admitting an inf-residuation) with $\OLR = \R\cup\cb{\pm\infty}$ using the 'inf-addition' $\isum$ (see \cite{HamelSchrage12, RockafellarWets98}). The inf-residuation on $\OLR$ is given by
\[
r \idif s  = \inf\cb{t\in\R \mid r \leq s \isum t}
\]
for all $r,s\in\OLR$, compare \cite{HamelSchrage12} for further details.
\end{example}

Basic notions from real analysis can be easily extended to set--valued functions mapping into the conlinear space $\G^\triup$, for instance a function $f \colon X \to \G^\triup$ is called convex when
\begin{equation*}
\label{EqConvFct} \forall x_1, x_2 \in X, \; \forall t \in \of{0,1} \colon
 f\of{tx_1+(1-t)x_2} \lel tf\of{x_1} \isum \of{1-t}f\of{x_2}.
\end{equation*}
Moreover $f$ is called positively homogeneous when
\[
\forall 0<t, \forall x \in X \colon f\of{tx} \lel tf\of{x},
\]
and it is called sublinear if it is positively homogeneous and convex. As a standard notation, we refer to the image set of a subset $A\subseteq X$ through $f$ by 
$$f\sqb{A}=\cb{f(x)\in \G^\triup\st x\in A}\subseteq \G^\triup$$
 and the (effective) domain of a function $f:X\to \G^\triup$ is the set 
$$\dom f=\cb{x\in X\st f(x)\neq \sup \G^\triup}.$$
 A function $f:X\to \G^\triup$ is called proper if $\dom f\neq\emptyset$ and $\inf \G^\triup\notin f\sqb{X}$.

%When dealing with set--valued functions, scalarization is a common tool for optimization problems. 
We recall that the recession cone of a nonempty closed convex set $A\subseteq Z$ is the closed convex cone $0^+A=\cb{z\in Z\st A+\cb{z}\subseteq A}$,
compare \cite[p.6]{Zalinescu02}. By definition, $0^+\emptyset=\emptyset$ is assumed. If $A\in\G^\triup\setminus\cb{\emptyset}$, then $0^+A=A\idif A$ and $C\subseteq 0^+A$ are satisfied. Especially, $\Int\of{0^+A}\neq\emptyset$ and $(0^+A)^-\subseteq C^-$, hence $B^*\cap (0^+A)^-$ is a weak$^*$ compact base of $(0^+A)^-$.

Each element of $\G^\triup$ is closed and convex and satisfies $A=A+C$, hence by a separation argument we can prove
\begin{align}\label{eq:scal_representation_Set}
\forall A\in \G^\triup:\quad A=\bigcap\limits_{z^*\in B^*}\cb{z\in Z\st -\sigma(z^*| A)\leq -z^*(z)},
\end{align}
where $\sigma(z^*| A)=\sup\cb{z^*(z)\st z\in A}$ is the support function of $A$ at $z^*$. Therefore, $A=\emptyset$ if and only if there exists a $z^*\in B^*$ such that $-\sigma(z^*| A)=+\infty$, or equivalently if the same holds true for all $z^*\in B^*$.

According to this notation, introducing the family of scalarizations for $f \colon X \to \G^\triup$ as the extended real-valued functions $\vp_{f, z^*} \colon  X \to \R\cup\cb{\pm\infty}$ defined by
\[
\forall z^* \in \C:\quad  \vp_{f,z^*}\of{x} = \inf\cb{-z^*\of{z} \mid z \in f\of{x}}
\]
we obtain from \eqref{eq:scal_representation_Set} the following representation of $f$
\[
\forall x \in X: \quad f\of{x} = \bigcap_{z^* \in B^*}\cb{z \in Z \st \vp_{f, z^*}\of{x} \leq -z^*\of{z}}.
\]
Some properties of $f$ are inherited by its scalarizations and vice versa. For instance, $f$ is convex if and only if  $\vp_{f, z^*}$ is convex for each $z^* \in B^*$.

To some extent, continuity or its relaxations are a common assumption in variational inequality applications to optimization. The following definition summarizes those continuity concepts that are used in the sequel.

\begin{definition}\label{def:l.s.c.}
\begin{enumerate}[(a)]
\item
Let $\vp:X\to\OLR$ be a function, $x_0\in X$. Then $\vp$ is said to be lower semicontinuous (l.s.c.) at $x_0$ iff
\begin{equation*}
\forall r\in\R:\quad r<\vp_{f,z^*}(x_0)\;\Rightarrow\;
\exists U\in\mathcal U_X:\, \forall u\in U:\,  r<\vp_{f,z^*}(x_0+u).
\end{equation*}
\item
A set $\Psi=\cb{\vp_i:X\to\OLR\st i\in I}$ is lower equicontinuous in $x_0\in \bigcap\limits_{i\in I}\dom\vp_i$ if
\[
\forall \eps>0\,\exists U\in\mathcal U_X\,\forall x\in x_0+U\,\forall i\in I:\quad
\vp_i(x_0)\leq \vp_i(x)+\eps
\]
\item
Let $\psi:S\subseteq X\to Z$ be a function, then $\psi$ is called $C$--upper continuous ($C$--continuous) at $x_0\in S$ iff
\[
\forall V\in \mathcal U_Z\;\exists U\in\mathcal U_X\;\forall x\in S\cap (U+x_0):\quad \psi(x)\in \psi(x_0)+V+C;
\]
\item
Let $F:X\to \P(Z)$ be a function, then $F$ is called upper Hausdorff continuous at $x_0\in S$ iff
\[
\forall V\in \mathcal U_Z\;\exists U\in\mathcal U_X\;\forall x\in  x_0+U:\quad F(x)\subseteq F(x_0)+V;
\]
\item
Let $f:X\to\G^\triup$ be a function, $M^*\subseteq\C$. Then $f$ is said $M^*$-- lower semicontiuous ($M^*$--l.s.c.) at $x_0$ iff
$\vp_{f,z^*}$ is l.s.c. at $x_0$ for all $z^*\in M^*$.
\item
Let $f:X\to\G^\triup$ be a function. If
\[
f(x)\lel \liminf\limits_{u\to 0}f(x+u)=\bigcap\limits_{U\in\mathcal U}\cl\co \bigcup\limits_{u\in U}f(x+u)
\]
is satisfied, then $f$ is lattice lower semicontinuous (lattice l.s.c.) at $x$.
\end{enumerate}
\end{definition}

In \cite{HeydeSchrage11R}, it has been proven that if $f$ is $\C$--l.s.c. at $x$, then it is also lattice l.s.c. at $x$.
As we assume $\Int C\neq\emptyset$, $ B^*$ exists and $f$ is $\C$--l.s.c. at $x$ if and only if $f$ is $ B^*$--l.s.c. at $x$.
One can show that if $f$ is convex, then $f$ is lattice l.s.c.
if and only if $\gr f=\cb{(x,z)\st z\in f(x)} \subseteq X \times Z$ is a closed set with respect to the product
topology, see \cite{HamelSchrage13}.

In \cite[Proposition 2.4]{CrespiRoccaSchrage14W}  it has been proven that
$f:X\to\G^\triup$ is upper Hausdorff continuous at $x_0\in\dom F$ if and only if $\Psi=\cb{\vp_{f,z^*}:X\to\OLR\st z^*\in B^*}$ is lower equicontinuous at $x_0$ in which case $f$ is $B^*$- lower semicontinuous at $x_0$ which in turn implies lower lattice continuity at $x_0$,
 compare  \cite{HeydeSchrage11R} for a detailed study of continuity concepts for set--valued functions.

\begin{remark}\label{rem:Epi_ext}
In this paper we mainly refer to $\G^\triup$--valued functions. However this is not a restriction as any set--valued function $F: X\to \P(Z)$ can be associated with its  $\G^\triup$--valued extension given by $F^C:X\to\G^\triup$ defined by
$$
F^C(x)=\left\{\begin{array}{lcl}
\cl\co\of{F(x)+C}, & \mbox{\rm if} & F(x)\neq\emptyset\\
F(x)=\emptyset & & \mbox{\rm elsewhere}.
\end{array}\right.
$$
Recalling that $C$--convexity of $F$ is defined by
\[
\forall x,y\in X,\;\forall t\in\of{0,1}:\quad tF(x)+(1-t)F(y)\subseteq F(tx+(1-t)y)+C.
\]
we have that $F^C$ is convex if $F$ is $C$--convex.
\end{remark}

In order to apply the set--valued results to vector--valued functions, we need the following definition of set--valued extension of a vector--valued function.
\begin{definition}
Let $\psi:S\subseteq X\to Z$ be a vector--valued function. The $\G^\triup$ extension of $\psi$ is the set--valued function $\psi^C:X\to\G^\triup$, 
\[
\forall x\in X:\quad \psi^C(x)=\left\{\begin{array}{lcl}
\cb{\psi(x)}+C, & \mbox{\rm if} & x\in\dom\psi\\
\emptyset & & \mbox{\rm otherwise}.
\end{array}\right.
\]
\end{definition}
Obviously, $\dom \psi^C=S$ and for all $z^*\in B^*$ it holds
$$
\vp_{\psi^C,z^*}(x)=\left\{\begin{array}{ll}
										-z^*\of{\psi(x)}\in\R & \mbox{\rm if \;}x\in S\\
										+\infty & \mbox{\rm elsewhere}.
													\end{array}
										\right.
$$

\begin{lemma}\label{lem:Equicont_CCont}
A vector--valued function $\psi:S\subseteq X\to Z$ is 
%For all $z^*\in B^*$ it holds 
%$$
%\vp_{\psi^C,z^*}(x)=\left\{\begin{array}{ll}
%										-z^*\psi(x)\in\R & \mbox{\rm if \;}x\in S\\
%										+\infty & \mbox{\rm elsewhere}.
%													\end{array}
%										\right.
%$$
 $C$--continuous at $x_0\in S$ if and only if $\Psi=\cb{\vp_{\psi^C,z^*}:X\to\OLR\st z^*\in B^*}$ is lower equicontinuous at $x_0$.
\end{lemma}
\proof
Assume first that $\psi$ is $C$--continuous at $x_0$, so 
\[
\forall V\in \mathcal U_Z\;\exists U\in \mathcal U_X\;\forall x\in \cb{x_0}+(U\cap S):\quad \psi(x)\in \cb{\psi(x_0)}+V+C.
\]
As $B^*$ is weak$^*$-compact, for each $\varepsilon>0$ there exists a $z_\varepsilon\in -\Int C$ such that $-\varepsilon=\inf\limits_{z^*\in B^*}-z^*(z_\varepsilon)$. Setting $V_\varepsilon=z_\varepsilon+C$, then $C$--continuity implies
\[
\forall \varepsilon>0\;\exists U\in \mathcal U_X\;\forall x\in \cb{x_0}+(U\cap S):\quad \psi(x)\in \psi(x_0))+V_\varepsilon+C,
\] 
or equivalently
\[
\forall \varepsilon>0\;\exists U\in \mathcal U_X\;\forall x\in \cb{x_0}+(U\cap S):\quad -z^*(\psi(x)-\psi(x_0))>-\varepsilon.
\] 
%as $\varphi_{(f,z^*)}(x)=-z^*(\psi(x))$. 
This implies lower equicontinuity of $\Psi$.

Now assume $\Psi$ is lower equicontinuous at $x_0$. By contradiction, let $\psi$ not be $C$--continuous at $x_0$. Then there exists a $V_0\in\mathcal U_Z$ and a net $\cb{x}_{i\in I}\subseteq S$ with $x_i\to x_0$,  such that $\psi(x_i)\notin \psi(x_0)+V_0+C$ for all $i\in I$.

As $B^*$ is weak${}^*$-compact, for all $V\in \mathcal U(0)$ it holds
\[
\sup\limits_{z^*\in B^*}\inf\limits_{z\in V}-z^*(z)<0,
\]
compare \cite[Remark 3.32]{HeydeSchrage11R}.
Especially, by a separation argument 
\[
0>-\varepsilon_0=\sup\limits_{z^*\in B^*}\inf\limits_{z\in V_0}-z^*(z)\in \R
\]
is true and, again by a separation argument, for each $i\in I$ there exists a $z^*_i\in B^*$ such that 
\[
 -z^*_i(\psi(x_0))+\inf\limits_{z\in V_0}-z^*_i(z)\geq -z^*_i(\psi(x_i))\in\R.
\] 
Thus,
\[
\forall i\in I:\quad -z^*_i(\psi(x_i))+\frac{1}{2}\varepsilon_0\lneq  -z^*_i(\psi(x_0)),
\] 
contradicting the lower equicontinuity of $\Psi$.
\pend

If additionally the ordering cone $C$ (and hence $C^{-}$) is polyhedral, then $\psi$ is $C$--continuous if and only if $\psi^C$ is $B^*$-l.s.c. at $x$ (see \cite[Corollary 5.6]{Luc89}). %Under the same assumption on $C$, lower semicontinuity of a finite set of scalarizations, namely those associated with the extreme directions of $\C$, is known to be equivalent to $C$--continuity of $\psi$.

Finally we define the restriction of a set--valued function $f:X\to\G^\triup$ to a segment with end points $x_0,x\in X$ as
$f_{x_0,x}:\R\to\G^\triup$, given by
\[
f_{x_0,x}(t)=\begin{cases}
	f(x_0+t(x-x_0)),\text{ if } t\in\sqb{0,1};\\
	\emptyset,\text{ elsewhere.}
\end{cases}
\]
%This is equivalent to the restriction of a scalar valued function $\vp:X\to\OLR$ to the same segment, defined by
%\[
%\vp_{x_0,x}(t)=\begin{cases}
%	\vp(x_t),\text{ if } t\in\sqb{0,1};\\
%	+\infty,\text{ elsewhere.}
%\end{cases}
%\]
Setting $x_t=x_0+t(x-x_0)$ for all $t\in\R$,
the scalarization $\vp_{f_{x_0,x},z^*}:\R\to\OLR$ of the restricted function $f_{x_0,x}$ is equal to the restriction $(\vp_{f,z^*})_{x_0,x}:\R\to\OLR$ of the scalarization of $f$ for all $z^*\in\C$.

%If $f$ is convex, $x_0,x_t\in \dom f$ for some $t\in\of{0,1}$, then $\of{\vp_{f,z^*}}_{x_0,x}$ is lower semicontinuous on $\of{0,t}$, hence $f_{x_0,x}$ is lattice l.s.c. on $\of{0,t}$.
An immediate generalization of the results in the remainder of this note is to replace convexity of $f$ by radial convexity of $f$ at $x_0$, meaning that $f_{x_0,x}:\R\to\G^\triup$ is convex for all $x\in X$ and likewise replacing lower semicontinuity by the corresponding radial definition.
In \cite{CrespiHamelSchrage13W} and \cite{CrespiRoccaSchrage14W}, the convexity assumption is dropped and replaced by more general monotonicity assumptions on the scalarization of the set--valued function.

\section{Minimality and variational inequality formulation}

In set optimization several notions of minimimality can be defined through the order introduced in $\G^\triup$. In this paper we focus on the following definitions that introduce two different notions, a stronger and a weaker one, respectively.

\begin{definition}\label{def:minimizer} \cite{HeydeLoehne11}
Let $f:X\to\G^\triup$ be a function. Then $x_0\in\dom f$ is called a {\em minimizer} of $f$ if the following holds true.
\begin{align}\label{eq:Min}
\tag{$Min$}
\forall x\in X:\quad \of{f(x)\lel f(x_0)\;\Rightarrow f(x)=f(x_0)}.
\end{align}
\end{definition}

This definition corresponds to the so-called set criterion which became
popular due to the work of Kuroiwa and collaborators, compare e.g. \cite{Kuroiwa1998natural}.

\begin{definition}\label{def:weak_minimizer}
Let $f:X\to\G^\triup$ be a function. Then $x_0\in\dom f$ is called a {\em weak l}, {\em scalarized weak} or {\em weak  minimizer} of $f$ if either $f(x)=Z$, or
\begin{align}\label{eq:w-l-Min}
\tag{$w$-$l$-$Min$}
\forall x\in X:\quad f(x_0)\nsubseteq \Int f(x);
\end{align}
\begin{align}\label{eq:w-sc-Min}
\tag{$w$-$sc$-$Min$}
\forall x\in X\,\exists z^*\in B^*:\quad \vp_{f,z^*}(x_0)\leq \vp_{f,z^*}(x)\neq-\infty;
\end{align}
\begin{align}\label{eq:w-Min}
\tag{$w$-$Min$}
\forall x\in X\,\forall U\in\mathcal U:\quad f(x_0)\oplus U\nsubseteq f(x).
\end{align}
\end{definition}
The chain of implications in Definition \ref{def:weak_minimizer} is \eqref{eq:w-l-Min}$\Rightarrow$\eqref{eq:w-sc-Min}$\Rightarrow$\eqref{eq:w-Min} (see \cite[Proposition 2.9]{CrespiSchrage13b}) hence each weak-l-minimizer of $f$ in the sense of \cite{HernandezMarin05} is a weak minimizer. The $l$ in the definition of weak l minimizers refers to the specific ordering in use, sometimes called the 'lower' set ordering, in contrast to the upper ordering, compare also \cite{Kuroiwa1998natural}.
If $f=F^C:X\to\G^\triup$ and $F(x_0)$ is a compact set, then the three types of weak minimizers coincide \cite[Proposition 2.1]{CrespiRoccaSchrage14W}.
The motivation of our naming lays in the special case $f=\psi^C$. Then $x_0\in\dom f$ is a weak minimizer of $f$ if and only if $\psi(x_0)$ is a weakly efficient element of $\psi\sqb{X}$, i.e. for all $x\in\dom f$ it holds $\psi(x_0)\notin\psi(x)+\Int C$. Likewise, $x_0$ is a minimizer of $f$ if and only if $\psi(x_0)$ is an efficient element of $\psi\sqb{X}$, i.e.  for all $x\in\dom f$, $\psi(x_0)\in\psi(x)+ C$ implies $\psi(x)\in\psi(x_0)+ C$.

To introduce a variational inequality associated with the set optimization of $f:X\to\G^\triup$, we first need a notion of derivative of $f$. Recent results on scalar and vector Minty type variational inequalities such as \cite{crespiJOTA, CreGinRoc08} have used the concept of (lower) Dini derivative to state the problem. The structure of a $\inf$--residuated image space allows to propose such a derivative also for set--valued maps. Notably, the definition extends the Dini derivative of scalar--valued functions to extended real--valued functions (see e.g. \cite{HamelSchrage13, Diss}).

\begin{definition}\label{def:DirDer}
Let $Y$ be a $\inf$--residuated order complete conlinear space, $f:X\to Y$ and $x,u\in X$. The {\em upper} and {\em lower Dini directional derivative} of $f$ at $x$ in direction $u$ are given by
\begin{align*}
f^\uparrow(x,u)&=\limsup\limits_{t\downarrow 0}\frac{1}{t}\of{f(x+tu)\idif f(x)}
=\inf\limits_{0<s}\sup\limits_{0<t\leq s}\frac{1}{t}\of{f(x+tu)\idif f(x)};\\
f^\downarrow(x,u)&=\liminf\limits_{t\downarrow 0}\frac{1}{t}\of{f(x+tu)\idif f(x)}
=\sup\limits_{0<s}\inf\limits_{0<t\leq s}\frac{1}{t}\of{f(x+tu)\idif f(x)}.
\end{align*}
If both derivatives coincide, then $f'(x,u)=f^\uparrow(x,u)=f^\downarrow(x,u)$ is the {\em Dini directional derivative} of $f$ at $x$ in direction $u$.
\end{definition}

\begin{remark}\label{ex:scalar_DirDer}
Definition \ref{def:DirDer} actually provides a generalization of the classical notion of Dini derivative for scalar--valued functions. Indeed let $\vp:X\to\OLR$ be an extended real--valued scalar function. If $\vp(x+tu)\in\R$ is satisfied for all $t\in\sqb{0,t_0}$ for a given $0<t_0$, then the differential quotient is real, too, hence in this case the above defined derivatives coincide with the standard definition in the literature, compare \cite{giorgi1992dini}.
If $x\notin\dom \vp$, then $\vp(x+tu)\idif \vp(x)=-\infty$ for all $t>0$, so $\vp'(x,u)=-\infty$.
On the other hand if $\vp(x)=-\infty$, then $\vp(x+tu)\idif \vp(x)=-\infty$ whenever $\vp(x+tu)=-\infty$ and $\vp(x+tu)\idif \vp(x)=+\infty$ else. The value of the derivatives in this case depends on the behaviour of $\vp$ in the proximity of $x$.
\end{remark}

As the Dini derivatives are defined 'radially', it is easy to see that the following statements hold under radial assumptions, too.
For notational simplicity, we refrain from this generalization, hoping to improve the clarity of the general scheme presented.

The following characterization of the Dini derivative extends a classical result to set-valued functions.

\begin{proposition}\label{prop:DirDer_of_conv}\cite[Proposition 3.4]{CrespiSchrage13b}
Let $Y$ be a $\inf$--residuated order complete conlinear space, $f:X\to Y$. If $f$ is convex, then the Dini derivative exists for all $x,u\in X$ and it holds
\[
f'(x,u)=\inf\limits_{0<t}\frac{1}{t}\of{f(x+tu)\idif f(x)}.
\]
Moreover, $f':X\times X\to Y$ is sublinear in its second component.
\end{proposition}

If $Y=\G^\triup$, then for all $x,u\in X$ and $0<s$ the directional derivative of a convex function $f:X\to\G^\triup$ is
\[
f'(x,u)=\cl\bigcup\limits_{0<t\leq s}\frac{1}{t}\of{f(x+tu)\idif f(x)},
\]
the differential quotient is decreasing as $t$ converges towards $0$. Moreover, as $\Int C\neq\emptyset$ is assumed,
\[
\Int f'(x,u)=\bigcup\limits_{0<t\leq s}\Int \frac{1}{t}\of{f(x+tu)\idif f(x)}
\]
is satisfied for all $x,u\in X$ and all $0<s$, compare \cite[Lemma 3.5]{CrespiSchrage13b}.

\begin{proposition}
Let $f:X\to\G^\triup$ be a convex, set--valued function, then
\[
\liminf\limits_{t\downarrow 0}\frac{1}{t}\of{f(x+tu)\idif f(x)}=\bigcap\limits_{0<s}\cl\bigcup\limits_{0<t<s}\frac{1}{t}\of{f(x+tu)\idif f(x)},
\]
the upper Painleve Kuratowski limit of the difference quotient, compare \cite[p. 21]{Aubin71}.
\end{proposition}
\proof
Indeed, we only need to check the convexity of the set $\bigcup\limits_{0<t<s}\frac{1}{t}\of{f(x+tu)\idif f(x)}$.
Let $z_1, z_2\in \bigcup\limits_{0<t<s}\frac{1}{t}\of{f(x+tu)\idif f(x)}$ be given, then there exists $0<t_1,t_2<s$ such that
\[
f(x)+t_iz_i\subseteq f(x+t_iu)
\]
is true for $i=1,2$. Let $r\in\sqb{0,1}$ be given, $t_0=(1-r)t_1+rt_2$. By convexity of the set $f(x)$ it holds
\[
f(x)+t_0\of{\frac{(1-r)t_1}{t_0}z_1+\frac{rt_2}{t_0}z_2}=(1-r)\of{f(x)+t_1 z_1}+r\of{f(x)+t_2z_2},
\]
hence by convexity of the function $f$
\[
f(x)+t_0\of{\frac{(1-r)t_1}{t_0}z_1+\frac{rt_2}{t_0}z_2}\subseteq (1-r)f(x+t_1u)+rf(x+t_2u)\subseteq f(x+t_0u),
\]
implying $\co\cb{z_1,z_2}\subseteq \bigcup\limits_{0<t<s}\frac{1}{t}\of{f(x+tu)\idif f(x)}$.

\pend

We will frequently make use of the following relations between the Dini derivative of a set--valued function and those of its scalarization.

\begin{proposition}\label{prop:vp'fLeqf'}\cite[Proposition 2.36]{CrespiSchrage13a}
Let $f:X\to\G^\triup$ be a convex function, $x,u\in X$. Then
\begin{align*}
&\bigcap\limits_{z^*\in B}\cb{z\in Z\st \vp'_{f,z^*}(x,u)\leq -z^*(z)}\lel f'(x,u);\\
&\forall z^*\in \C:\quad \vp'_{f,z^*}(x,u)\leq -\sigma(z^*|f'(x,u)) .
\end{align*}
\end{proposition}

Although in general the scalarization of the derivative is not equal to the derivative of the scalarization, the equality occurs in the special case of the epigraphical extension of a vector--valued function.

\begin{proposition}\label{prop:reg_sharp}\cite[Proposition 3.10]{CrespiSchrage13b}
Let $\psi:S\subseteq X\to Z$ be a $C$--convex vector--valued function, $f=\psi^C:X\to\G^\triup$ its epigraphical extension, $x,u\in X$.
Then for all $z^*\in\C$ it holds
\begin{equation}\label{eq:reg_sharp}\tag{SR}
\forall z^*\in \C:\quad -\sigma(z^*|f'(x,u))=\vp'_{f,z^*}(x,u).
\end{equation}
\end{proposition}

For a general function $f:X\to \G^\triup$, if \eqref{eq:reg_sharp} is satisfied, then also the weaker condition
\begin{equation}\label{eq:reg_weak}\tag{WR}
f'(x,u)= \bigcap\limits_{z^*\in B^*}\cb{z\in Z\st \vp'_{f,z^*}(x,u)\leq -z^*(z)}
\end{equation}
holds true. Notably, when $f:X\to\G^\triup$ is the epigraphical extension of a $C$--convex vector function $\psi:S\subseteq X\to Z$, then Property \eqref{eq:reg_weak} is satisfied.

In the sequel, property \eqref{eq:reg_sharp} will be referred to as strong regularity, while property \eqref{eq:reg_weak} will be referred to as weak regularity.

Attempts to characterize minimizers and weak minimizers in set optimization through variational inequalities have been made in \cite{CrespiHamelSchrage13W, CrespiSchrage13a} and \cite{CrespiRoccaSchrage14W, CrespiSchrage13b}.
We recall the definitions of Stampacchia and Minty variational inequalities and their scalarizations used in the previous papers.

\begin{definition}\label{def:set_Stamp}\cite[Definition 3.12]{CrespiSchrage13a}
Let $f:X\to\G^\triup$ be a convex function, $x_0\in\dom f$. Then $x_0$ solves the set--valued Stampacchia inequality iff
\begin{equation}\label{eq:set_Stamp}
\tag{$SVI_M$}
f(x_0)=Z\;\vee\;\forall x\in \dom f:\, f(x)\neq f(x_0)\,\Rightarrow\, 0\notin f'(x_0,x-x_0).
\end{equation}
\end{definition}

\begin{definition}\label{def:scalar_Stamp}\cite[Definition 3.14]{CrespiSchrage13a}
Let $f:X\to\G^\triup$ be a convex function, $x_0\in\dom f$. Then $x_0$ solves the scalarized Stampacchia inequality iff
\begin{equation}\label{eq:scalar_Stamp}
\tag{$svi_M$}
f(x_0)=Z\;\vee\;\forall x\in \dom f:\, f(x)\neq f(x_0)\,\Rightarrow\,\exists z^*\in B^*:\, 0< \vp'_{f,z^*}(x_0,x-x_0).
\end{equation}
\end{definition}

\begin{definition}\label{def:set_Minty}\cite[Definition 3.18]{CrespiSchrage13a}
Let $f:X\to\G^\triup$ be a convex function, $x_0\in\dom f$. Then $x_0$ solves the set--valued Minty inequality iff
\begin{equation}\label{eq:set_Minty}
\tag{$MVI_M$}
\forall x\in X:\;f(x)\neq f(x_0)\quad\Rightarrow\quad  f'(x,x_0-x)\not\subseteq 0^+f(x).
\end{equation}
\end{definition}

\begin{definition}\label{def:scalar_Minty}\cite[Definition 3.19]{CrespiSchrage13a}
Let $f:X\to\G^\triup$ be a convex function, $x_0\in\dom f$. Then $x_0$ solves the scalarized Minty inequality iff
\begin{equation}\label{eq:scalar_Minty}
\tag{$mvi_M$}
\forall x\in X:\;f(x)\neq f(x_0)\quad\Rightarrow\quad\exists z^*\in B^*:\quad \vp_{f,z^*}(x)\neq-\infty\,\wedge\, \vp'_{f,z^*}(x,x_0-x)< 0.
\end{equation}
\end{definition}

In \cite{CrespiSchrage13a} the following scheme has been proved for convex set--valued functions.

\begin{proposition}\label{prop:image1}
 Let $f:X\to\G^\triup$ be a convex function, $x_0\in\dom f$.
\begin{enumerate}[(a)]
\item
The following implications hold without further assumptions.
\[
\eqref{eq:scalar_Stamp}\;\Rightarrow\; \eqref{eq:set_Stamp}\;\Rightarrow\;\eqref{eq:Min}\;\Rightarrow\;
\eqref{eq:scalar_Minty}\;\Leftarrow\; \eqref{eq:set_Minty};
\]
\item
If the weak regularity assumption \eqref{eq:reg_weak} is satisfied, then $\of{\eqref{eq:scalar_Stamp}\;\Leftrightarrow\; \eqref{eq:set_Stamp}}$ is true, strong regularity \eqref{eq:reg_sharp} implies $\of{\eqref{eq:scalar_Minty}\;\Leftrightarrow\; \eqref{eq:set_Minty}}$
\item
If $f$ is $B^*$--l.s.c. in $x_0$ and there exists a finite subset $M^*$ of $B^*$ such that for all $x\in X$ the defining inequality in \eqref{eq:scalar_Minty} is attained at some element $z^*\in M^*$, then $\of{\eqref{eq:Min}\Leftrightarrow\eqref{eq:scalar_Minty}}$ is true.
\item
Especially if $\psi:S\subseteq X\to Z$ is given, $f(x)=\psi^C(x)$ for all $x\in X$, then
\[
\eqref{eq:scalar_Stamp}\;\Leftrightarrow\; \eqref{eq:set_Stamp}\;\Rightarrow\;\eqref{eq:Min}\;\Rightarrow\;
\eqref{eq:scalar_Minty}\;\Leftrightarrow\; \eqref{eq:set_Minty}
\]
is satisfied. If additionally the ordering cone $C$ is polyhedral and $\psi$ is $C$--continuous at $x_0$, then the following scheme is true.
\[
\eqref{eq:scalar_Stamp}\;\Leftrightarrow\; \eqref{eq:set_Stamp}\;\Rightarrow\;\eqref{eq:Min}\;\Leftrightarrow\;
\eqref{eq:scalar_Minty}\;\Leftrightarrow\; \eqref{eq:set_Minty}.
\]

\end{enumerate}
\end{proposition}

The assumption of a finite set $M^*\subseteq B^*$ satisfying \eqref{eq:scalar_Minty} to prove \eqref{eq:Min} cannot be dropped, as the following example shows.

\begin{example}
Let $X=\R$ and $Z=l^\infty$ be given with the usual ordering cone $C=\cb{z\in Z\st \forall n\in \N\, z_n\geq 0}$. The function $f=\psi^C:X\to\G(Z,C)$ is defined by $\dom \psi=\sqb{-1,1}$ and
\[
\forall x\in\dom \psi\,\forall n\in N:\quad
(\psi(x))_n=
\max\cb{(\sqrt{n^2-1}-n)(x+1),\frac{1}{n-\sqrt{n^2-1}}(x-1)}.
\]
Then $-e_n^*=(0,...,0,-1,0,...)\in \C$ and $\vp_{f,-e_n^*}(x)=(\psi(x))_n$ is true for all $n\in \N$ and all $x\in\sqb{-1,1}$ and $f$ is convex. % and radially upper Hausdorff continuous in $1$. 
However,
\[
\forall -1< x<1:\quad f(1)\subsetneq f(x),
\]
i.e. $1$ is not a minimizer of $f$, while
\[
\vp'_{f,-e_n^*}(x,1)=
\begin{cases}
\sqrt{n^2-1}-n,&\text{ if } x<\frac{\sqrt{n^2-1}}{n};\\
\frac{1}{n-\sqrt{n^2-1}},&\text{ if } x\geq\frac{\sqrt{n^2-1}}{n}.
\end{cases}
\]
As directional derivatives are positively homogeneous, this implies \eqref{eq:scalar_Minty} is satisfied at $1$.
\end{example}

It is left as an open question whether or not other assumptions would be sufficient to provide equivalence between \eqref{eq:scalar_Minty} and \eqref{eq:Min}.

\medskip
The implications provided in Proposition \ref{prop:image1} are illustrated in the following figure.
   \begin{center}
     \begin{figure}[h!]
        \centering
        \includegraphics[width=14cm]{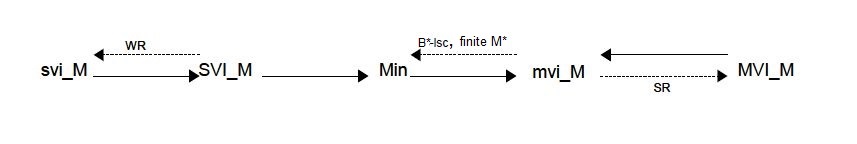}
    \end{figure}
    \end{center}

Weaker inequalities can be introduced as well to characterize weak efficiency.

\begin{definition}\label{def:weak_set_Stamp}
Let $f:X\to\G^\triup$ be a convex function, then $x_0$ solves the weak set--valued Stampacchia inequality iff
\begin{equation}\label{eq:weak_set_Stamp}
\tag{$SVI_W$}
f(x_0)=Z\;\vee\;\forall x\in X:\, 0\notin \Int f'(x_0,x-x_0).
\end{equation}
\end{definition}

\begin{definition}\label{def:weak_scalar_Stamp}
Let $f:X\to\G^\triup$ be a convex function, $x_0\in\dom f$. Then $x_0$ solves the weak scalarized Stampacchia inequality iff
\begin{equation}\label{eq:weak_scalar_Stamp}
\tag{$svi_W$}
f(x_0)=Z\;\vee\;\forall x\in X:\,\exists z^*\in B^*:\,
0\leq \vp'_{f,z^*}(x_0,x-x_0).
\end{equation}
\end{definition}

\begin{definition}\label{def:weak_set_Minty}
Let $f:X\to\G^\triup$ be a convex function, $x_0\in\dom f$. Then $x_0$ solves the weak set--valued Minty inequality iff
\begin{equation}\label{eq:weak_set_Minty}
\tag{$MVI_W$}
f(x_0)=Z\;\vee\;\forall x\in X:\, f'(x,x_0-x)\nsubseteq \Int 0^+f(x).
\end{equation}
\end{definition}

\begin{definition}\label{def:weak_scalar_Minty}
Let $f:X\to\G^\triup$ be a convex function, $x_0\in\dom f$. Then $x_0$ solves the weak scalarized Minty inequality, iff
\begin{equation}\label{eq:weak_scalar_Minty}
\tag{$mvi_W$}
f(x_0)=Z\;\vee\;\forall x\in X:\,\exists z^*\in B^*:\quad \vp_{f,z^*}(x)\neq-\infty\,\wedge\, \vp'_{f,z^*}(x,x_0-x)\leq 0.
\end{equation}
\end{definition}

In \cite{CrespiRoccaSchrage14W, CrespiSchrage13b}, the following scheme has been proved for convex set--valued functions $f:X\to\G^\triup$.

\begin{proposition}\label{prop:image2}
 Let $f:X\to\G^\triup$ be a convex function, $x_0\in\dom f$.
\begin{enumerate}[(a)]
\item
The following implications hold without further assumptions.
\[
\eqref{eq:weak_scalar_Stamp}\;\Rightarrow\; \eqref{eq:w-sc-Min}\;\Rightarrow\;
\eqref{eq:weak_scalar_Minty}\;\Leftarrow\; \eqref{eq:set_Minty};
\]
\[
\eqref{eq:weak_scalar_Stamp}\;\Rightarrow\; \eqref{eq:weak_set_Stamp}\;\Leftrightarrow\;\eqref{eq:w-Min};
\]
\item
If the strong regularity assumption \eqref{eq:reg_sharp} is satisfied, then the following implications hold with equivalence.
\[
\eqref{eq:weak_set_Stamp}\;\Leftrightarrow\;\eqref{eq:weak_scalar_Stamp}\;\Leftrightarrow\; \eqref{eq:w-sc-Min} %\;\Leftrightarrow\; \eqref{eq:w-Min}
\]
\item
If the set $B^*$ in \eqref{eq:w-sc-Min} can be replaced by a finite subset $M^*\subseteq B^*$,  then the following equivalence is satisfied.
\[
\eqref{eq:weak_scalar_Stamp}\Leftrightarrow\eqref{eq:w-sc-Min}
\]
\item
If the set $B^*$ in \eqref{eq:weak_scalar_Minty} can be replaced by a finite subset $M^*\subseteq B^*$ and $f$ is $M^*$--l.s.c. in $x_0$, then
\[
\eqref{eq:weak_scalar_Stamp}\;\Leftrightarrow\;\eqref{eq:w-sc-Min}\;\Leftrightarrow\;\eqref{eq:weak_scalar_Minty}
\]
\item
If $f=F^C$ with $F(x_0)\subseteq Z$ compact, then
\[
\eqref{eq:weak_set_Stamp}\;\Leftrightarrow\;\eqref{eq:weak_scalar_Stamp}\;\Leftrightarrow\;\eqref{eq:w-sc-Min}\of{\;\Leftrightarrow\;\eqref{eq:w-Min}\;\Leftrightarrow\;\eqref{eq:w-l-Min}}.
\]
\item
If $f=F^C$ with $F(x_0)\subseteq Z$ compact, the scalarizations $\vp^\triup_{f,z^*}$ are proper for all $z^*\in B^*$ and $f_{x_0,x}$ is upper Hausdorff continuous for all $x\in\dom f$, then
\[
\eqref{eq:weak_set_Stamp}\;\Leftrightarrow\;\eqref{eq:weak_scalar_Stamp}\;\Leftrightarrow\;\eqref{eq:w-Min}\;\Leftrightarrow\;\eqref{eq:weak_scalar_Minty}.
\]
\item
If $Z$ is of finite dimension or $C$ is Daniell, then if $\psi:S\subseteq X\to Z$ is given, $f(x)=\psi^C(x)$ for all $x\in X$, then
\[
\eqref{eq:weak_set_Stamp}\;\Leftrightarrow\; \eqref{eq:weak_scalar_Stamp}
\;\Leftrightarrow\;\eqref{eq:w-Min}\;\Leftrightarrow\;\eqref{eq:w-sc-Min}
\;\Rightarrow\;\eqref{eq:scalar_Minty}\;\Leftrightarrow\; \eqref{eq:set_Minty}
\]
is satisfied. If additionally  $\psi_{x_0,x}$ is $C$-continuous for all $x\in S$, %and the ordering cone $C$ is polyhedral, 
then all implications are satisfied with equivalence.
\begin{align*}
\eqref{eq:set_Stamp}\;\Leftrightarrow\; \eqref{eq:scalar_Stamp}\;\Leftrightarrow\;\eqref{eq:w-Min}\;&\Leftrightarrow\;\eqref{eq:w-sc-Min}\\ &\Leftrightarrow\;
\eqref{eq:w-l-Min}\;\Leftrightarrow\;\eqref{eq:scalar_Minty}\;\Leftrightarrow\; \eqref{eq:set_Minty}.
\end{align*}

\end{enumerate}
\end{proposition}

In general a solution to the scalarized Minty variational inequality \eqref{eq:scalar_Minty} is not a solution to the set-valued one, as the following example shows.

\begin{example}
Let $X=\R$, $Z=\R^2$ and $C=\R^2_+$ be given, $x_0=\frac{2}{3}$ and $f:\R\to\G^\triup$ with
\[
f(x)=
\begin{cases}
\cb{z=(z_1,z_2)^T\in\R^2\st z_1+z_2\geq (1-\frac{1}{2}x),\, z_1\geq x,\, z_2\geq x }, &\text{ if } 0\leq x\leq \frac{2}{3};\\
\emptyset,& \text{ otherwise.}
\end{cases}
\]
Then $f(0)=\cb{(x,1-x)^T\st 0\leq x\leq 1}+C$ is the sum of a compact set and $C$, $f$ is convex and upper Hausdorff continuous in the domain and each scalarization is proper.
Let $z^*=(-1,-1)^T$, then $\vp^\triup_{f,z^*}(0)=1$ and $\of{\vp^\triup_{f,z^*}}'(0,1)=-\frac{1}{2}$ while $f'(0,1)=\of{1,1}^T+C\subseteq \Int 0^+f(0)$. So \eqref{eq:scalar_Minty} (and \eqref{eq:weak_scalar_Minty}) is satisfied even with a finite subset of $B^*$ but \eqref{eq:set_Minty} (and \ref{eq:weak_set_Minty}) is not satisfied.
Especially, \eqref{eq:reg_sharp} is not satisfied.
\end{example}

%{\color{red}
%\begin{example}
%SVI does not imply svi
%\end{example}
%}
\medskip
The implications provided in Proposition \ref{prop:image2} are illustrated in the following figure.
   \begin{center}
     \begin{figure}[h]
        \centering
        \includegraphics[width=14cm]{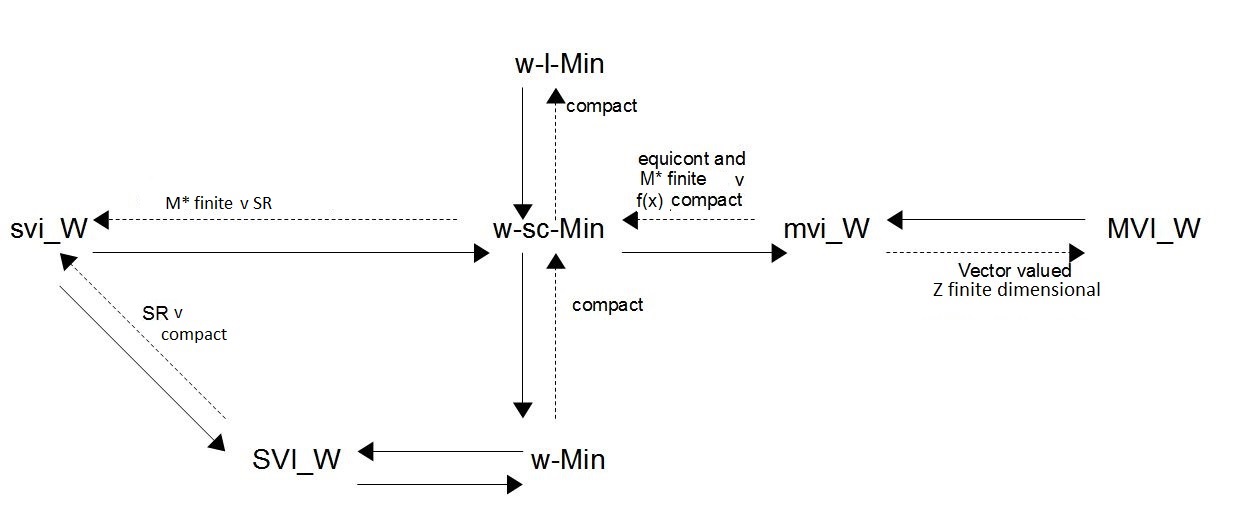}
    \end{figure}
    \end{center}

\section{Main Results}

While minimizers clearly are weak minimizers, we still need to prove that the same implication holds between the strong and the weak formulation of the variational inequalities. This result closes the loop between the previous schemes, providing a complete wreath between variational inequalities and minimality.

The following results prove the relations holding between the four couples of inequalities introduced in the previous section.

\begin{proposition}\label{prop:4_1}
Let $f:X\to\G^\triup$ be a convex function, $x_0\in\dom f$. If $x_0$ solves \eqref{eq:set_Stamp}, then it also solves \eqref{eq:weak_set_Stamp} and $f(x)=f(x_0)$ implies that $f_{x_0,x}$ is constant on $\sqb{0,1}$.
\end{proposition}
\proof
Let $x_0$ solve \eqref{eq:set_Stamp}. If $x\in\dom f$ and $f(x)\neq f(x_0)$, then by definition $0\notin f'(x_0,x-x_0)$, implying $0\notin\Int f'(x_0,x-x_0)$, as required in \eqref{eq:weak_set_Stamp}.
Assume $f(x)=f(x_0)$ and $f(x_t)\neq f(x)$ for some $t\in\of{0,1}$. By convexity, $f(x_0)\subsetneq f(x_t)$ is satisfied, hence
\[
0\in f(x)\idif f(x_0)\subseteq f'(x_0,x-x_0).
\]
On the other hand, the derivative is positively homogeneous, hence
\[
(1-t)f'(x_0,x-x_0)=f'(x_0,x_t-x_0)
\]
and thus by \eqref{eq:set_Stamp}
\[
0\notin f'(x_0,x_t-x_0),
\]
a contradiction. Hence in this case $f_{x_0,x}$ is constant on $\sqb{0,1}$ and
\begin{eqnarray*}
f'(x_0,x-x_0)=0^+f(x_0).
\end{eqnarray*}
In this case, $0\in\Int f'(x_0,x-x_0)$ implies $f(x_0)=Z$, the set--valued weak Stampacchia inequality is satisfied for all $x\in \dom f$.

If $x\notin \dom f$, then either $\dom f_{x_0,x}\cap\of{0,1}=\emptyset$ and
$f'(x_0,x-x_0)=\emptyset$,
or there exists a $t\in\of{0,1}$ such that $f(x_t)\neq \emptyset$. In this case, the same argument as above proves the statement, replacing $x$ by $x_t$.
\pend

\begin{proposition}
Let $f:X\to\G^\triup$ be a convex function, $x_0\in\dom f$. If $x_0$ solves \eqref{eq:scalar_Stamp}, then it also solves \eqref{eq:weak_scalar_Stamp} and $f(x)=f(x_0)$ implies that $f_{x_0,x}$ is constant on $\sqb{0,1}$.
\end{proposition}
\proof
Assume $f(x)=f(x_0)$ and $f(x_t)\neq f(x)$ for some $t\in\of{0,1}$. By convexity, $f(x_0)\subsetneq f(x_t)$ is satisfied, hence
\[
(1-t)\vp'_{f,z^*}(x_0,x-x_0)=\vp'_{f,z^*}(x_0,x_t-x_0)\leq \vp_{f,z^*}(x_t)\idif \vp_{f,z^*}(x_0)\leq 0
\]
is satisfied for all $z^*\in B^*$ and \eqref{eq:scalar_Stamp} implies the existence of $\bar z^*\in  B^*$ such that
\[
 0<\vp'_{f,\bar z^*}(x_0,x_t-x_0).
\]
But this implies $\vp_{f,\bar z^*}(x_0)<\vp_{f,\bar z^*}(x_t)$, a contradiction.
Hence in this case $f_{x_0,x}$ is constant on $\sqb{0,1}$ and
\begin{eqnarray*}
\forall z^*\in B^*:\quad \vp_{f,z^*}(x_0)=-\infty\,\vee\, \vp'_{f,z^*}(x_0,x-x_0)=0.
\end{eqnarray*}
In this case, $\vp_{f,z^*}(x_0)=-\infty$ for all $z^*\in B^*$ implies $f(x_0)=Z$, the scalarized weak Stampacchia inequality is satisfied for all $x\in\dom f$.

If $x\notin \dom f$, then either $\dom f_{x_0,x}\cap\of{0,1}=\emptyset$ and
$\vp'_{f,z^*}(x_0,x-x_0)=+\infty$
for all $z^*\in  B^*$, or there exists a $t\in\of{0,1}$ such that $f(x_t)\neq \emptyset$. In this case, the same argument as above proves the statement, replacing $x$ by $x_t$.
\pend

\begin{proposition}
Let $f:X\to\G^\triup$ be a convex function, $x_0\in\dom f$. If $x_0$ solves \eqref{eq:scalar_Minty}, then it solves \eqref{eq:weak_scalar_Minty}.
\end{proposition}
\proof
Under the assumption of  \eqref{eq:scalar_Minty}, let $f(x)=f(x_0)$ be satisfied. By convexity, $\vp_{f,z^*}(x_t)\leq \vp_{f,z^*}(x)$ is true for all $z^*\in B^*$ and all $t\in\of{0,1}$. Thus either $f_{x_0,x}$ is constant on $\sqb{0,1}$, in which case $f(x_0)=Z$ or $\vp'_{f,z^*}(x,x_0-x)=0$ for all $z^*\in B^*$ with $\vp_{f,z^*}(x_0)\neq -\infty$, or there exists a $t\in\of{0,1}$ such that $f(x_t)\supsetneq f(x)$. In this case, by assumption there exists a $z^*\in B^*$ such that $-\infty\neq \vp_{f,z^*}(x_t)\leq \vp_{f,z^*}(x)$ and $\vp'_{f,z^*}(x_t,x_0-x_t)<0$.
By convexity of $f$ this implies $\vp'_{f,z^*}(x,x_0-x)<0$.
\pend

\begin{proposition}\label{prop:4_4}
Let $f:X\to\G^\triup$ be a convex function, $x_0\in\dom f$. If $x_0$ solves \eqref{eq:set_Minty}, then it solves \eqref{eq:weak_set_Minty}.
\end{proposition}
\proof
Under the assumption of  \eqref{eq:set_Minty}, let $f(x)=f(x_0)$ be satisfied. By convexity, $f(x_t)\lel f(x)$ is true for all $t\in\of{0,1}$. Thus either $f_{x_0,x}$ is constant on $\sqb{0,1}$, in which case $f(x_0)=Z$ or $f'(x,x_0-x)=0^+f(x)\nsubseteq \Int 0^+f(x)$, or there exists a $t\in\of{0,1}$ such that $f(x_t)\supsetneq f(x)$. In this case, by assumption  $f'(x_t,x_0-x_t)\nsubseteq 0^+f(x_t)$. Let $s\in \of{0,1}$, then
\begin{align*}
f'(x,x_0-x)&\supseteq \frac{1}{s+t-st}\of{f(x_t+s(x_0-x_t))\idif f(x)}\\
			  &\supseteq \frac{1}{s+t-st}\of{\of{f(x_t+s(x_0-x_t))\idif f(x_t)}\oplus \of{f(x_t)\idif f(x)}}\\
%			 &\supseteq\of{1-\frac{t}{s+t-st}}\frac{1}{s}\of{f(x_t+s(x_0-x_t))\idif f(x_t)}\oplus \frac{t}{s+t-st}\frac{1}{t}\of{f(x_t)\idif f(x)}.
\end{align*}
By assumption, $f(x)\subseteq f(x_t)$, hence
\[
0^+f(x_t)\subseteq \of{f(x_t)\idif f(x)},
\]
which implies
\begin{align*}
&\of{f(x_t+s(x_0-x_t))\idif f(x_t)}\oplus \of{f(x_t)\idif f(x)}\\
\supseteq
&\of{f(x_t+s(x_0-x_t))\idif f(x_t)}\oplus 0^+f(x_t)\\
=
&f(x_t+s(x_0-x_t))\idif f(x_t)
\end{align*}
and therefore
\[
f'(x,x_0-x)\supseteq \frac{s}{s+t-st}\of{\frac{1}{s} f(x_t+s(x_0-x_t))\idif f(x_t)}.
\]

Moreover, $f'(x_t,x_0-x_t)\nsubseteq 0^+f(x_t)$, hence choosing $s\in\of{0,1}$ small enough,
\[
\frac{1}{s}\of{f(x_t+s(x_0-x_t))\idif f(x_t)}\nsubseteq 0^+f(x)
\]
is satisfied, proving
\begin{align*}
f'(x,x_0-x)
\nsubseteq \frac{s}{s+t-st} 0^+f(x),
\end{align*}
for some $s\in\of{0,1}$. Thus $f'(x,x_0-x)\nsubseteq \Int 0^+f(x)$, as desired.
\pend

%\pagebreak

The following easy example may serve to illustrate that the implications  proven in Propositions \ref{prop:4_1} to Proposition \ref{prop:4_4} may not be reversed.

\begin{example}
Let $X=Z=\R^2$ be ordered by the usual ordering cone $C=R^2_+$ and $S\subseteq X$ given as
\[
S=\R^2_+\cap\cb{x\in X\st x_1+x_2\geq 0}.
\]
The function $\psi:S\to Z$ is given as $\psi(x)=x$, hence the set--valued extension is $\psi^C:X\to\P(Z)$ with $\psi^C(x)=x+C$ whenever $x\in\dom \psi^C=S$ and $\psi^C(x)=\emptyset$, elsewhere.
The set of weak minimizers is the boundary of $S$, while the set of minimizers is the set $\cb{x\in S\st x_1+x_2=1}$:
 Let $x_0=(0,2)^T$ be fixed, then it is an easy task to prove that $x_0$ does not satisfy \eqref{eq:set_Minty} and hence does not solve \eqref{eq:scalar_Minty}, \eqref{eq:scalar_Stamp}, \eqref{eq:set_Stamp} or \eqref{eq:Min}.
 On the other hand, as $x_0$ is a weak l minimizer, it satisfies all the weak versions of variational inequalities presented in the course of this paper.
\end{example}

\medskip
 Merging the results of Propositions \ref{prop:4_1} to \ref{prop:4_4} with those presented in Proposition \ref{prop:image1} and Proposition \ref{prop:image2}, the following figure illustrates the obtained implications, thus summarizing the result of the presented paper.

   \begin{center}
     \begin{figure}[h]
        \centering
        \includegraphics[width=14cm]{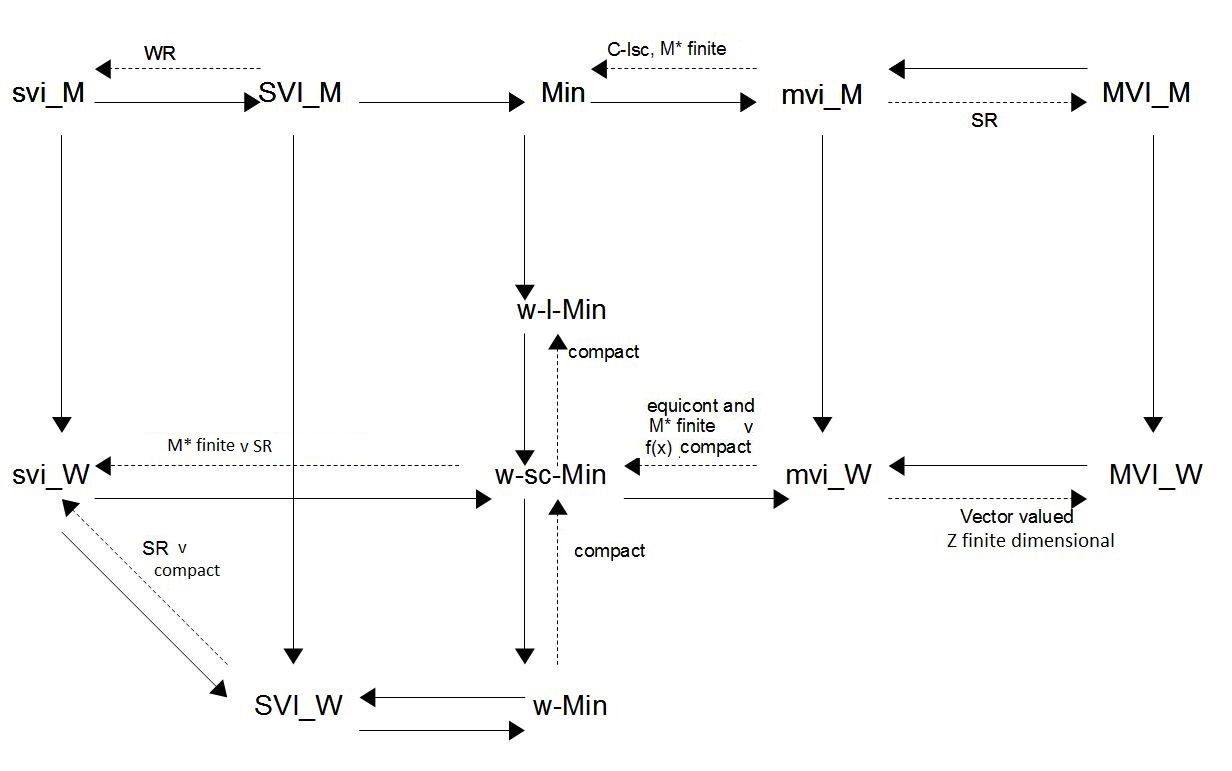}
    \end{figure}
    \end{center}

%If additionally$f=\psi^C$ is assumed, i.e. $f$ is the epigraphical extension of a vector--valued function, then the solution sets of scalarized and set--valued versions of each variational inequality coincide.

  \pagebreak
   
 The following scheme represents the implications proven for vector--valued functions and their set--valued extensions.
   \begin{center}
     \begin{figure}[!h]
        \centering
        \includegraphics[width=14cm]{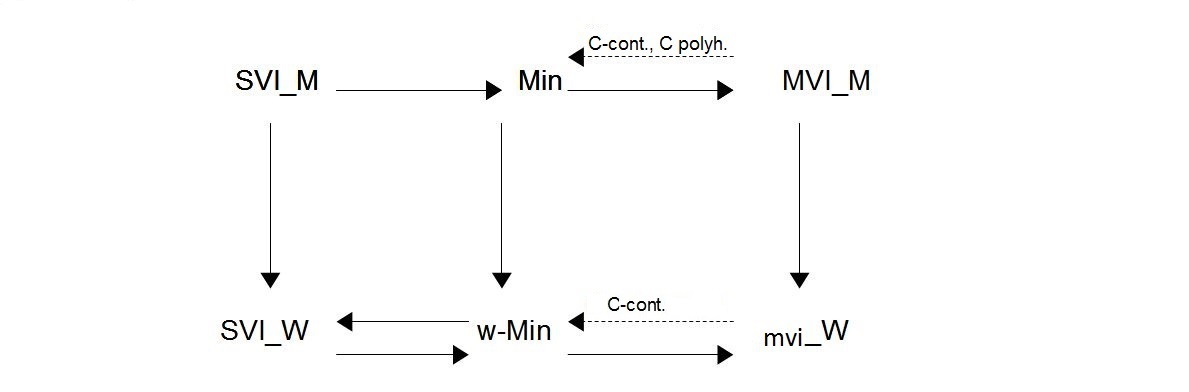}
    \end{figure}
    \end{center}

    In case $Z$ is of finite dimension or $C$ is Daniell, then the scalarized weak Minty variational inequality and the weak Minty variational inequality coincide, too.
		
We remark that the latter scheme is a straightforward extension of the scheme of relations originally provided by Giannessi for his vector variational inequalities in vector optimization. Therefore, as an application, we have proved that set optimization approach provides a useful tool to study vector optimization, by considering the epigraphical extension of the objective function. The main advantage we see in this approach is to work in an order complete space, rather than a partially ordered space.

\section*{Acknowledgement}

The authors are extremely thankful to the editor and two anonymous referees who have provided several suggestions to improve the paper to its current form.

\addcontentsline{toc}{section}{Bibliography}
%\bibliography{../../../../Literatur}

\end{document}